\documentclass[12pt]{article}

\usepackage{algorithmic}

\usepackage{amsthm,amsmath,amssymb}

\usepackage{stmaryrd}

\usepackage[colorlinks=true,citecolor=black,linkcolor=black,urlcolor=blue]{hyperref}

\usepackage{pstricks}
\newcommand{\showgrid}{}
\newcommand{\gridon}{\renewcommand{\showgrid}{\psset{subgriddiv=1,griddots=10,gridlabels=6pt}\psgrid}}
\gridon 

\usepackage{graphicx}

\usepackage{color}
\usepackage{multido}

\theoremstyle{plain}
\newtheorem{thm}{Theorem}

\theoremstyle{definition}
\newtheorem{dfn}[thm]{Definition}
\newtheorem{ex}[thm]{Example}
\newtheorem{con}[thm]{Conjecture}

\theoremstyle{remark}

\newgray{hellgrau}{.89}

\def\bit{\begin{itemize}}
\def\eit{\end{itemize}}
\def\beq{\begin{equation}}
\def\eeq{\end{equation}}
\def\bpf{\begin{proof}}
\def\epf{\end{proof}}

\def\EM#1{{\em #1}}

\def\figref#1{Figure~\ref{#1}}
\def\secref#1{Section~\ref{#1}}

\def\Z{{\mathbb Z}}

\def\1{{\mathbf 1}}
\def\0{{\mathbf 0}}

\def\defeq{:=}

\def\of#1{\left(#1\right)}
\def\absof#1{\left|#1\right|}
\def\setof#1{\left\{#1\right\}}

\def\pas#1{\left(#1\right)}

\def\symm{\mathfrak S}

\def\inv{{\mathbf{inv}}}
\def\ASM{alternating sign matrix}
\def\ASMes{alternating sign matrices}
\def\DPP{descending plane partition}

\usepackage[applemac]{inputenc}

\title{\bf A bijection between permutation matrices and \DPP s
  without special parts, which respects the quadruplet of statistics
  considered by Behrend, Di Francesco and Zinn--Justin.}

\author{Markus Fulmek\thanks{Research supported by the National Research Network ``Analytic
Combinatorics and Probabilistic Number Theory'', funded by the
Austrian Science Foundation.}\\
\small Faculty of Mathematics\\[-0.8ex]
\small University of Vienna\\[-0.8ex] 
\small Vienna, Austria\\
\small\tt Markus.Fulmek@univie.ac.at\\
}

\date{{2018--06--07}}

\def\asmset#1{\mathcal{A}_{#1}}
\def\dppset#1{\mathcal{D}_{#1}}
\def\asmsubset#1{\mathcal{A}^0_{#1}}
\def\dppsubset#1{\mathcal{D}^0_{#1}}
\def\stati{i} 
\def\statp{p} 
\def\stats{s} 
\def\statb{q} 
\def\quadruplet{\pas{\statp, \stati, \stats, \statb}}

\def\Iverson#1{\left[#1\right]}
\def\numtext#1{$\#\pas{\text{{\small #1}}}$}

\def\ppcell{{rd--cell}}

\def\mmcell{{lu--cell}}

\def\secA#1{\section{#1}}
\def\secB#1{\subsection{#1}}
\def\secC#1{\subsubsection{#1}}

\def\makroPlus{%
\pscircle[linecolor=black,linewidth=0.025](0,0){0.25}%
\psline[linecolor=black,linewidth=0.075](-0.18,0)(0.18,0)%
\psline[linecolor=black,linewidth=0.075](0,-0.18)(0,0.18)%
}
\def\makroMinus{%
\pscircle[linecolor=red,linewidth=0.025](0,0){0.25}%
\psline[linecolor=red,linewidth=0.075](-0.2,0)(0.2,0)%
}
\begin{document}

\maketitle


\begin{abstract}
  We present a bijection between permutation matrices and \DPP s
  without special parts, which respects the quadruple of statistics
  considered by Behrend, Di Francesco and Zinn--Justin. This bijection
  involves the inversion words of permutations 
  and an representation of \DPP s as families of non--intersec\-ting lattice paths.
\end{abstract}

\secA{Introduction}
\label{sec:intro}

It is a well--known fact \cite{Zeilberger:1995,Kuperberg:1996,Bressoud:1999} that the enumeration
\bit
\item
of \DPP s with parts not exceeding $n$ (let us denote the set of these objects by $\dppset{n}$)
\item and of \ASMes\ of dimension $n$ (let us denote the set of these objects by $\asmset{n}$)
\eit
gives the same number:
\begin{equation}
\label{eq:nof-ASMs}
\absof{\asmset{n}} = \absof{\dppset{n}} = \prod_{j=0}^{n-1}\frac{\pas{3\cdot n +1}!}{\pas{n + j}!}.
\end{equation}
\secB{Search for a bijection}
It appears to be quite difficult to find some ``natural'' bijection
$$\Phi:\mathcal{A}_n\to\mathcal{D}_n.$$
However, there are two additional informations which might help
in the search for such bijection: 
\bit
\item There is a \EM{quadruplet} $\quadruplet$ of statistics for both $\asmset{n}$ and $\dppset{n}$,
	i.e. there are functions
	$$
	\sigma_A:\asmset{n}\to\Z^4 \text{ and } \sigma_D:\dppset{n}\to\Z^4,
	$$
	such that for all $y\defeq\quadruplet\in\Z^4$ the preimages are \EM{equinumerous}, i.e.,
	$$
	\absof{\sigma_A^{-1}\of y} = \absof{\sigma_D^{-1}\of y}
	$$
	(see \cite[Theorem 1]{Behrend:2013}, the details
	are given in the next section).
\item There are certain subsets of $\asmset{n}$ and of $\dppset{n}$, namely
	\bit
	\item \ASMes\ with statistic $\stats=0$ (let us denote this set by $\asmsubset{n}$;
		it is, in fact, the set of $n\times n$ \EM{permutation matrices}), 
	\item \DPP s with statistic $\stats=0$ (let us denote this set by $\dppsubset{n}$),
	\eit
	which are \EM{much simpler} to understand and for which it is, in fact, quite easy to give ``natural''
	bijections (see below).
\eit
So one obvious approach would be to search for a bijection $\Phi$ which \EM{respects} this quadruplet
of statistics $\quadruplet$; in the sense that for all $A\in\asmset{n}$ there should hold:
{\small
\begin{equation}
\label{eq:respect-quadruplet}
\sigma_A\of{A} = \sigma_D\of{\Phi\of A}\text{ for all }A\in\asmset{n}.
\end{equation}
}

Clearly, such bijection $\Phi$ restricted to the subset $\asmsubset{n}\subset\asmset{n}$ would give a bijection
$$\Psi:\asmsubset{n}\to \dppsubset{n}.$$

So if we find such ``restricted'' bijection
$\Psi$ which respects the quadruplet of statistics in the sense of \eqref{eq:respect-quadruplet},
then we might hope to ``extend'' it somehow to the desired ``full'' bijection $\Phi$.

The purpose of this note is to present a simple bijection $\Psi$ which indeed
respects the quadruplet of statistics in the sense of \eqref{eq:respect-quadruplet}:
The construction of $\Psi$ relies
on the representation of \DPP s as families of non--intersecting lattice
paths and on a certain ``visualization'' of the statistic $\stati$ (as number of certain
cells in an \ASM).

\secB{Other bijections}
It should be noted that there are other bijections $\asmsubset{n}\to\dppsubset{n}$:
Maybe the simplest one was already mentioned by Lalonde \cite[p.981]{Lalonde:2006},
who referred to the inversion table of a permutation: This table (also called \EM{inversion word}) gives a (well--known) 
unique encoding of permutations. The
missing link from inversion words to $\dppsubset{n}$ was explained by Striker \cite[Lemma 5]{Striker:2011}.
Striker uses  monotone triangles as intermediate objects to establish the bijection between
inversion words and $\asmsubset{n}$ (we shall call this Striker's bijection), but this intermediate
step is not necessary: Instead, we can employ directly
the well--known
encoding of permutations by inversion words (we shall call this Lalonde's bijection). 
Unfortunately, none of these two simple bijections respects the statistic $\statb$ (see \figref{fig:non-respecting}).
Ayyer \cite{Ayyer:2010}  presented another (inductively constructed) bijection,
which does not respect the statistic $\stati$ (see \cite[p. 1786]{Ayyer:2010}).

\secB{Organization of this note}
This note is organized as follows:
\bit
\item \secref{sec:background} contains basic definitions and
background information,
\item \secref{sec:dpp-rep} presents a representation of \DPP s as families of
non--intersecting lattice paths (in the ``obvious'' way: for the expert it will
suffice to look at \figref{fig:lalonde-example}),
\item \secref{sec:asm-inv} presents a ``visualization'' of inversions in \ASMes\ 
	used for our bijection,
\item \secref{sec:nonspecial-bijection} presents a bijection $\Psi:\asmsubset{n}\to \dppsubset{n}$
	 which respects the quadruplet $\quadruplet$ of statistics.
\eit

\secA{Background information}
\label{sec:background}

For reader's convenience, we recall some background information.
\subsection{Descending plane partitions}
\label{sec:background-dpp}
Here is the definition of \DPP s as given by Mills, Robbins and Rumsey \cite[Definitions 2--4]{MillsRobbinsRumsey:1983}:
\begin{dfn}[\DPP]
\label{dfn:dpp}
A \EM{\DPP} is an array $\pi=\pas{a_{i,j}}$, $1\leq i\leq j<\infty$, of
\EM{positive integers} 
$$
\pi=\;\;
\begin{matrix}
a_{1,1} & a_{1,2} & a_{1,3} &         & \cdots       &        & \cdots      & a_{1,\mu_1} \\
        & a_{2,2} & a_{2,3} &         & \cdots       & \cdots & a_{2,\mu_2} & \\
        &         &         &         & \cdots       &        &             & \\
        &         &         &         & \cdots       &        &             & \\
        &         &         & a_{k,k} & \cdots &  a_{k,\mu_k} &             &
\end{matrix}
$$
such that
\begin{enumerate}
\item rows are \EM{weakly decreasing}, i.e., $a_{i,j}\geq a_{i,j+1}$ for all $i=1,\dots,k$ and $i\leq j < \mu_i$,
\item columns are \EM{strictly decreasing}, i.e., $a_{i,j}> a_{i+1,j}$ for all $i=1,\dots,k-1$ and $i<j\leq \mu_{i+1}$,
\item $a_{i,i}>\mu_i-i+1$ for all $i=1,\dots,k$,
\item $a_{i,i}\leq \mu_{i-1}-i+2$ for all $i=2,\dots,k$.
\end{enumerate}

Clearly, conditions 3 and 4 imply
$$
\mu_1\geq\mu_2\geq\cdots\geq\mu_k\geq k.
$$
The \EM{parts} of a \DPP\ are the \EM{numbers} (with repetitions) that appear in the array.
The \EM{empty} array, which we denote by $\emptyset$, is explicitly allowed. 

A \DPP\ $\pi$ where no part is greater than $n$ (i.e., $\pi$ has at most $n-1$ rows)
is said to have dimension $n$. (So a \DPP\ of dimension $n$ may also be viewed as
a \DPP\ of dimension $k$, for all $k>n$.)

By the \EM{length} of row $i$ in \DPP\ $\pi$ we define the number of parts it contains (i.e. $\mu_i-i+1$).
So we may rephrase conditions 3 and 4 as follows:
\bit
\item[$3^\prime$.] The first part of row $i$ is greater than the length of row $i$ for $i=1,\dots,k$,
\item[$4^\prime$.] The first part of row $i$ is less or equal than the length of the preceding row ${i-1}$ for $i=2,\dots,k$.
\eit

A part $a_{i,j}$ in a \DPP\ is called \EM{special} if it does not exceed the number of parts to
its left (in its row $i$), i.e., if
$$
a_{i,j}\leq j-i.
$$
\end{dfn}

\begin{ex}
A typical example is the array
$$
\begin{matrix}
6 & 6 & 6 & 4 & {\red\underline 2} \\
  & 5 & 3 & {\red\underline 2} & {\red\underline 1} \\
  &   & 2 &   & 
\end{matrix}
$$
with $3$ rows and $10$ parts (written in descending order)
$$
6,6,6,5,4,3,2,{\red\underline 2,\underline 2,\underline 1},
$$
three of which are  special parts (indicated as underlined numbers;
note that the $2$ in the last row is \EM{not}
a special part):
$$
{\red \underline 2,\underline 2,\underline 1}.
$$
(This is the example $D_0$ considered by Lalonde \cite[Fig.~1]{Lalonde:2003}.)
\end{ex}

From now on, we shall use the shortcut \EM{DPP} for \DPP s.
\def\DPP{DPP}

\subsection{Alternating sign matrices}
Here is the definition of \ASMes\ as given by Mills, Robbins and Rumsey, see \cite[Definition 1]{MillsRobbinsRumsey:1983}:
\begin{dfn}[\ASM]
An \EM{\ASM} of dimension $n$ is an $n\times n$ square matrix which satisfies
\bit
\item all entries are $1$, $-1$ or $0$,
\item every row and column has sum $1$,
\item in every row and column the nonzero entries alternate in sign.
\eit

Suppose that $A=\pas{A_{i,j}}_1^n$ is an alternating sign matrix of dimension $n$.
Then the number of \EM{inversions} in $A$ is defined to be \cite[p.\ 344]{MillsRobbinsRumsey:1983}
\begin{equation}
\label{eq:inversions}
\sum_{\substack{1\leq i<k\leq n\\1\leq l<j\leq n}}A_{i,j}\cdot A_{k,l}.
\end{equation}

\end{dfn}

\begin{ex}
The following matrix is an example of an \ASM\ of dimension $5$:
$$
\begin{bmatrix}
0 & 1 & 0 & 0 & 0 \\
0 & 0 & 1 & 0 & 0 \\
1 & -1 & 0 & 1 & 0 \\
0 & 1 & 0 & -1 & 1 \\
0 & 0 & 0 & 1 & 0.
\end{bmatrix}
$$
\end{ex}
From now on, we shall use the shortcut \EM{ASM} for \ASMes.
\def\ASM{ASM}
\def\ASMes{ASMs}

\subsection{The Mills--Robbins--Rumsey conjecture}
Here is the Conjecture of Mills, Robbins and Rumsey \cite[Conjecture 3]{MillsRobbinsRumsey:1983},
slightly rephrased to fit our exposition:
\begin{con}
\label{con:MRR}
Suppose that $n, \statp, \stati, \stats$ are nonnegative integers, $0\leq \statp\leq n-1$.
Let $\asmset{n}\of{\statp, \stati, \stats}$ be the set of \ASMes\ such that
\begin{enumerate}
\item the size of the matrix is $n\times n$ (i.e., its dimension is $n$),
\item the number of $0$'s to the left of the $1$ in the first row is $\statp$,
\item the number of $-1$'s in the matrix is $\stats$,
\item the number of inversions in the matrix is $\stati+\stats$.
\end{enumerate}
On the other hand, let $\dppset{n}\of{\statp, \stati, \stats}$ be the set of \DPP s such that
\begin{enumerate}
\item no part exceeds $n$ (i.e., the dimension of the \DPP\ is $n$),
\item there are exactly $\statp$ parts equal to $n$,
\item there are exactly $\stati$ special parts,
\item there are a total of $\stati+\stats$ parts.
\end{enumerate}
Then $\asmset{n}\of{\statp, \stati, \stats}$ and $\asmset{n}\of{\statp, \stati, \stats}$
have the same cardinality --- the sets $\asmset{n}$ and $\dppset{n}$ are equidistributed
with respect to the \EM{triplet of statistics} $\pas{\statp, \stati, \stats}$.
\end{con}
This conjecture was proved by Behrend, Di Francesco and Zinn--Justin \cite[Theorem 1]{Behrend:2012}.

\subsection{The fourth statistic given by Behrend, Di Francesco and Zinn--Justin}
In \cite{Behrend:2013}, Behrend, Di Francesco and Zinn--Justin presented a fourth
statistic $\statb$ for \ASM s and \DPP s and showed that \ASM s and \DPP s are equidistributed
with respect to the \EM{quadruplet of statistics} $\quadruplet$ \cite[Theorem 1]{Behrend:2013}:
This statistic $\statb$ is
\bit
\item for $n$--dimensional \ASMes\ equal to the number of $0$'s to the right of the $1$ in the last row,
\item for $n$--dimensional \DPP s equal to the number of parts $n-1$ plus the number of rows of length $n-1$.
\eit

\secB{Permutation matrices and inversions}
Let $\sigma\in\symm_n$ be a permutation of the first $n$ natural numbers $\setof{1,2,\dots,n}$.
\secC{Inversions of a permutation}
Recall that an \EM{inversion} of $\sigma$ is a pair $\pas{i,j}$ such that
$i<j$ and $\sigma\of i > \sigma\of j$. 
For the number $\inv\of\sigma$ of \EM{all} inversions of
$\sigma$ we have $0\leq\inv\of\sigma\leq \frac{n\cdot\pas{n+1}}2$.

We may assign to $\sigma$ its \EM{inversion word} $\pas{a_1, a_2,\dots, a_{n-1}}$, where
$a_k$ is the number of inversions $\pas{i,j}$ with $\sigma\of j=k$, $k = 1,2,\dots, n-1$.
Clearly we have $0\leq a_k\leq n-k$ and $a_1+a_2+\cdots+a_{n-1} = \inv\of\sigma$.

Considering the \EM{permutation word}
$$
\pas{\sigma\of 1,\sigma\of 2,\dots,\sigma\of n},
$$
of $\sigma$, the inversion word's $k$--th entry $a_k$ is simply the number of elements
\EM{to the left} of $k$ (in the permutation word) which are \EM{greater} than $k$, and
it is easy to see that \EM{every} word $\pas{b_1,b_2,\dots,b_{n-1}}$ with $0\leq b_k\leq n-k$
determines a \EM{unique} permutation: Inversion words are, in this sense, just another
``encoding'' for permutations.

\secC{Permutation matrices}
A permutation $\sigma\in\symm_n$ can be represented by an $n\times n$--matrix $M$ with
entries
$$
M_{i,j} = \delta_{\sigma\of i,j}
$$
(where $\delta_{x,y}$ denotes Kronecker's delta: $\delta_{x,y} = 1$ if $x=y$,
$\delta_{x,y} = 0$ if $x\neq y$). We call this matrix the \EM{permutation matrix}
of $\sigma$: Clearly, it contains precisely one entry $1$ in every row and column.
\begin{ex}
Let $n=6$ and $\sigma\in\symm_6$ be the permutation with permutation word 
$$
\sigma=\pas{3 5 2 4 6 1}.
$$
The corresponding permutation matrix is
$$
\begin{bmatrix}
0 & 0 & 1 & 0 & 0 & 0 \\
0 & 0 & 0 & 0 & 1 & 0 \\
0 & 1 & 0 & 0 & 0 & 0 \\
0 & 0 & 0 & 1 & 0 & 0 \\
0 & 0 & 0 & 0 & 0 & 1 \\
1 & 0 & 0 & 0 & 0 & 0
\end{bmatrix},
$$
and the corresponding inversion word is
$$
\pas{2,3,1,1,1}.
$$
\end{ex}
Note that every \EM{permutation matrix} is an \ASM\ (which does not contain entries
$-1$), and that the definition of inversions \eqref{eq:inversions} for \ASMes\ is a generalization
of the number of inversions of a permutation (see also \secref{sec:asminvs}).

\secA{Representation of \DPP s as lattice paths}
\label{sec:dpp-rep}

If some row $i$ in a \DPP\ $\pi=\pas{a_{i,j}}$ 
is \EM{shorter} than $a_{i,i} - 1$, i.e.,
$$
\delta = \pas{a_{i,i} - 1} - \pas{\mu_i - i + 1 } > 0,
$$
then we pad this row with $\delta$ trailing zeroes (so the \EM{length} of some row is
the number of \EM{non--zero} parts in that row.)

\begin{figure}
\caption{Lattice--path--encoding for \DPP s.}
\label{fig:lalonde-example}
Consider the following two \DPP s of dimension $6$:
\begin{center}
$
\begin{matrix}
6 & 6 & 6 & 4 & {\red\underline 2}\\
  & 5 & 3 & {\red\underline 2} & {\red\underline 1}\\
  &   & 2 &  & \\
\end{matrix}
$
\hfil
$
\begin{matrix}
6 & 6 & 6 & 5 & 0\\
  & 4 & 3 & {\red\underline 1} \\
  &   & 2 &  & \\
\end{matrix}
$
\end{center}
The pictures below show the non--intersecting lattice paths corresponding to the above \DPP s: 
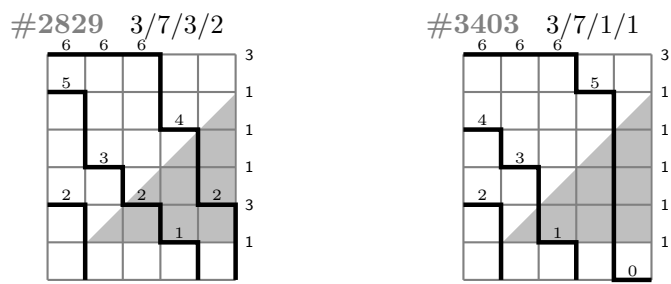
\begin{center}
\input ./save-graphics/DPP6ex1
\hfil
\input ./save-graphics/DPP6ex2
\end{center}
Note that the the right \DPP\ has a zero--padded first row, and that the 
special parts of these \DPP s correspond to the horizontal steps (at heights $>0$)
in the ``special range'' below the main diagonal $y=x$ (indicated by the gray triangle).
\end{figure}

Now we employ the well--known encoding of (shifted) tableaux as non--intersec\-ting
lattice paths; i.e., we encode a \DPP\
$\pi=\pas{a_{i,j}}$ of dimension $n$ with $r$ rows as an $r$--tuple of
non--intersecting lattice paths in the lattice $\Z^2$.
For reader's convenience, we shall describe the details of this encoding below, but the
idea can easily be obtained by looking at the illustrative
example in \figref{fig:lalonde-example}.

These lattice paths shall only use horizontal steps to the right or vertical steps downwards,
i.e., steps leading from lattice point $\pas{x,y}$ to lattice point $\pas{x+1,y}$ or to
lattice point $\pas{x,y-1}$.

The starting points of these $r$ lattice paths are the points
$$
S_i\defeq\pas{0, a_{i,i}},
$$
i.e., the lattice path corresponding to row $i$ starts 
on the vertical axis at height equal to the first part of row $i$.

The ending points of these lattice paths are the points 
$$
E_i\defeq\pas{a_{i,i}-1,0},
$$
i.e., the lattice path corresponding to row $i$ ends on the horizontal axis and
consists of $a_{i,i}-1$ horizontal steps at heights corresponding to the parts
of row $i$; including padded zeroes (if any), which correspond to steps at height $0$.

It is easy to see that the set $\dppset{n}$ of $n$--dimensional \DPP s is in bijection with
the set of nonintersecting lattice paths as defined above, with no starting point higher than
$n$. (But \EM{note} that in this representation, the number of horizontal steps of a path is equal to the length
of the corresponding row \EM{plus} the number of horizontal steps at height \EM{zero}.)

\secA{Inversions in \ASMes}
\label{sec:asm-inv}

\secB{Orientation of cells in \ASMes}

\begin{dfn}
By a \EM{cell} in some \ASM\ $A$ we simply mean the \EM{position} $\pas{i,j}$ (at
row $i$ and column $j$). If $A_{i,j}= 0$, we call this a \EM{zero--cell},
otherwise a \EM{non--zero cell}: So, non--zero cells can be either $1$--cells
or $\pas{-1}$--cells.

Observe that for every \EM{zero--cell} $\pas{i,j}$ in an \ASM\ there must be 
\bit
\item to the left or to the right (or both) of $\pas{i,j}$ a \EM{closest non--zero cell}
	in row $i$, precisely one of which must be a $1$--cell: If this $1$--cell lies to the
	left of $\pas{i,j}$, then call this cell \EM{left--oriented}, otherwise \EM{right-oriented}.
\item above or below (or both) $\pas{i,j}$ a \EM{closest non--zero} cell in column $j$,
	precisely one of which must be a $1$--cell: If this $1$--cell lies below
	$\pas{i,j}$, then call this cell \EM{down--oriented}, otherwise \EM{up-oriented}.
\eit
A zero--cell $\pas{i,j}$ is called
\bit
\item \EM{\ppcell}, if it is \EM{right--oriented} and \EM{down--oriented}.
\item \EM{\mmcell}, if it is \EM{left--oriented} and \EM{up--oriented}.
\eit
\end{dfn}

\secB{Inversions of \ASM s}
\label{sec:asminvs}
Observe that the \EM{\ppcell s} in the \EM{permutation matrix}
$A = \pas{\delta_{\sigma\of i,j}}$ 
of some permutation $\sigma\in\symm_n$
are in one--to--one correspondence with the \EM{inversions} of $\sigma$:
$\pas{i,j}$ is a \ppcell\ if 
\bit
\item entry $1$ in row $i$ is in column $y = \sigma\of i > j$ (to the right of column $j$),
\item and entry $1$ in column $j$ is in row $x = \sigma^{-1}\of j > i$ (below row $i$).
\eit
This is equivalent to
$$
\pas{\sigma\of i = y} > \pas{j = \sigma\of x}\text{ and } x > i,
$$
i.e., $\pas{i,x}$ is an \EM{inversion} of the permutation $\sigma$.

More generally, observe that the fourfold sum \eqref{eq:inversions}, which defines
the number of inversions of an \ASM, may
be rewritten as follows
$$
\sum_{\pas{i,l}=\pas{1,1}}^{\pas{n,n}} 
\sum_{\pas{k,j}=\pas{i+1,l+1}}^{\pas{n,n}} 
A_{i,j}\cdot A_{k,l}.
$$
Observe that the inner sum is simply the product
$$
\pas{\sum_{k=i+1}^{n}A_{i,j}}\cdot\pas{\sum_{j=l+1}^{n}A_{k,l}}.
$$
A moment's thought shows that this product is equal to
$$
\Iverson{\pas{i,j}\text{ is a $\pas{-1}$--cell or a \ppcell}}.
$$
Here, we used \EM{Iverson's notation}:
$$
\Iverson{\text{some assertion}}\defeq 
\begin{cases} 1 &\text{if the assertion is true,} \\ 0 & \text{ otherwise.}\end{cases}
$$
%
%
So the number of inversions of an \ASM\ $A$ (according to \eqref{eq:inversions})
is equal to the number of \ppcell s of $A$ \EM{plus} the number of
$\pas{-1}$--cells of $A$.

\secB{Quadruplet of statistics, reformulated}
These considerations immediately lead to the following reformulation of the
quadruplet of statistics $\quadruplet$:

\begin{center}
{\small
\begin{tabular}{|c|l|l|}
\hline
\multicolumn{3}{|c|}{Definition of the statistic $f$:} \\
\hline
$f$ & for $D\in\dppset{n}$: & for $A\in\asmset{n}$: \\
\hline
$\statp$ & \numtext{parts equal to $n$} & \numtext{\ppcell s in first row} \\
$\stati$ & \numtext{non--special parts} & \numtext{\ppcell s} \\
$\stats$ & \numtext{special parts} & \numtext{$\pas{-1}$--entries} \\
$\statb$ & \numtext{parts $(n-1)$} $+$ \numtext{rows of length $(n-1)$} & \numtext{\mmcell s in last row} \\
\hline
\end{tabular}
}
\end{center}
Now it is easy to show by a simple example that neither Lalonde's nor Striker's bijection
respect statistic $\statb$: Look at \figref{fig:non-respecting}, where we indicated
\bit
\item the $1$--cells by symbol
\psset{unit=0.5cm}\rput(0.3,0.2){\makroPlus}\hskip0.75em,
\item the $\pas{-1}$--cells by symbol
\psset{unit=0.5cm}\rput(0.3,0.2){\makroMinus}\hskip0.75em,
\item the \ppcell s by small equilateral right--angled triangles, where the sides of the
right angle point to the right and downwards,
\item and the \mmcell s by small equilateral right--angled triangles, where the sides of the
right angle point to the left and upwards.
\eit
\begin{figure}
\caption{Neither Lalonde's nor Striker's bijection respect statistic $\statb$.}
\label{fig:non-respecting}
Consider the $4$--dimensional \DPP\ $\pi=\begin{smallmatrix}4&4&4\\ &3&2\end{smallmatrix}$: $\statb\of\pi=2$,
since $\pi$ has one part $3$ and one row of length $3$.

The lattice path representation of $\pi$ is shown to
the left, the inversion word (needed for both Lalonde's and Striker's bijection) corresponding
to $\pi$ is $\pas{3,1,1}$, and the \EM{monotone triangle} (needed for Striker's bijection
\cite{Striker:2011}) corresponding to $\pi$ is shown on the top.

Both Striker's and Lalonde's bijection map $\pi$ to the \ASM\ shown to the right, with statistic $\statb$
equal to $3$ (the number of \mmcell s, indicated by south--east pointing triangles, in the last
row is $3$). However, $\pi$ should be
mapped to the \ASM\ shown at the bottom, which is the \EM{only} \ASM\ of dimension $4$ with quadruplet
of statistics $\pas{3,5,0,2}$.
\begin{center}
{\footnotesize
\begin{tabular}{ccccccc}
&&&{\bf 4}&&&\\
&&2&&{\bf 4}&&\\
&{\bf 2}&&{\bf 3}&&{\bf 4}&\\
1&&2&&3&&4\\
\end{tabular}
}
\end{center}
\begin{center}
\input ./save-graphics/DPPfig
\hfil
\input ./save-graphics/ASMfig1
\end{center}
\begin{center}
\input ./save-graphics/ASMfig2
\end{center}
\end{figure}

\secB{\ppcell s and \mmcell s are equinumerous}
Observe that for every \ASM\ $A$
the number of \ppcell s equals the number of \mmcell s.
We show this by a bijection (see \figref{fig:pp-mm-bijection}):
For some \ppcell\ $\pas{i,j}$, we construct two paths,
\bit
\item both starting at $\pas{i,j}$,
\item both proceeding only horizontally to the right or vertically downwards,
\item and both changing horizontal/vertical direction of movement whenever they encounter a non--zero cell.
\eit
One of these paths starts horizontally (to the right), the other starts vertically (downwards):
The pictures in \figref{fig:pp-mm-bijection} illustrate this simple idea. 
\begin{figure}
\caption{The bijection mapping \ppcell s to \mmcell s.}
\label{fig:pp-mm-bijection}
The left picture shows the two paths constructed for \ppcell\ $\pas{1,1}$: The first cell of crossing
of these paths is \mmcell\ $\pas{4,6}$, so $\pas{1,1}$ is mapped to $\pas{4,6}$.

The right picture shows the two paths constructed for \ppcell\ $\pas{3,2}$: The first cell of crossing
of these paths is \mmcell\ $\pas{6,5}$, so $\pas{3,2}$ is mapped to $\pas{6,5}$. Observe that these paths
also \EM{meet} (but do not \EM{cross}) at the $-1$--cell $\pas{5,4}$.
\begin{center}
\input ./save-graphics/ASM6bij1
\hfil
\input ./save-graphics/ASM6bij2
\end{center}
\end{figure}
Now observe that the paths starting horizontally must \EM{necessarily} end vertically, and vice versa.
Hence the paths \EM{must} eventually \EM{cross} at some cell, and any such cell of crossing must \EM{necessarily}
be a \mmcell\ (the paths could meet but not cross at some $\pas{-1}$--cell, see the right
picture in \figref{fig:pp-mm-bijection}): Map $\pas{i,j}$ to the \EM{first} cell of crossing $\pas{k,l}$
thus obtained. By symmetry (reflection at the second diagonal of $A$) it is immediately clear that this
construction gives a \EM{bijection}.


\secA{The bijection between \DPP s without special parts and permutation matrices}
\label{sec:nonspecial-bijection}

\secB{The statistic $\statb$ for \DPP s without special parts}
Observe that if a \DPP\ $\pi\in\dppset{n}$ \EM{without} special parts has a row of length $n-1$,
this row
\bit
\item must start with part $n$
\item and must not contain parts smaller than $n-1$
\eit
(this fact is immediately seen from the lattice path representation of $\pi$).

Stated otherwise: The somewhat complicated condition ``existence of a path of length $n-1$''
can be simply expressed as
\begin{equation}
\label{eq:path-non-special}
\text{\numtext{parts equal to $n$}} + \text{\numtext{parts equal to $n-1$}} \geq n-1.
\end{equation}
in the case of \DPP s without special parts.

\secB{The bijective construction}
Given some \EM{permutation matrix} $A$ (i.e., some \ASM\ without entries $-1$)
of dimension $n$, set $k=1$, $A_1=A$, and repeat the following step $n-1$ times:
\bit
\item Note down the number $a_k$ of \ppcell s in the first row of $A_k$,
\item delete the first row and the column containing the $1$ (i.e., column $a_k+1$),
\item rotate the matrix by $180^\circ$,
\item let $A_{k+1}$ be the $\pas{n-k}\times\pas{n-k}$--matrix thus obtained, and increase $k$ by $1$.
\eit
(See \figref{fig:bij-example} for an illustration of this construction.)

We claim that the sequence of $\pas{n-1}$ numbers $\pas{a_1, a_2,\dots, a_{n-1}}$ thus obtained
has the following properties:
\bit
\item $0\leq a_i\leq n-i$ for $i=1,2,\dots,n-1$, i.e., $\pas{a_1, a_2,\dots, a_{n-1}}$
	is an \EM{inversion word},
\item $a_1 = \statp\of A$,
\item $a_2 = \statb\of A - \Iverson{\statb\of A + \statp\of A \geq n}$\footnote{Here we use Iverson's notation again.},
\item $a_1+a_2+\cdots+a_{n-1} = \stati\of A$. 
\eit
Moreover, we claim that this construction gives a \EM{bijection} between permutation matrices and
inversion words.
\begin{proof}
It is clear that $0\leq a_k\leq n-k$, since $a_i$ is (by construction) the number
of \ppcell s in the first row of an $\pas{n-k+1}\times\pas{n-k+1}$ permutation matrix,
and this number cannot exceed $n-k$.

By construction, we also have $a_1=\statp\of A$.

Of the $\statb\of A$ \mmcell s in the last row of $A$, precisely \EM{one} will be deleted
in the first step of our construction \EM{if and only if} $\statb\of A + \statp\of A \geq n$: By the rotation
at the end of the first step, precisely these ``surviving'' \mmcell s will turn up as the
\ppcell s in the first row of the matrix at the beginning of the second step, whose number will give $a_2$.

Finally, if there are $a_k$ \ppcell s in the first row of the $\pas{n+1-k}\times\pas{n+1-k}$--matrix $A_k$
at the beginning of step $k$, then the column $a_k+1$ of $A_k$ contains precisely $a_k$ \mmcell s, since
the submatrix given by rows $2$ to $n+1-k$ and columns $1$ to $a_k$ must contain \EM{precisely} $a_k$ $1$--cells.
So in every step, the number of deleted \ppcell s \EM{equals} the number of deleted \mmcell s,
and after $n-1$ steps,
\bit
\item \EM{all} of the \ppcell s and \mmcell s of $A$ (recall that their number is \EM{twice} the
	number of \ppcell s) have been deleted,
\item and \EM{half} of these cells have been noted down during the construction; their number
	is equal to $a_1+\cdots+a_{n-1}$.
\eit
It is easy to see that the mapping from permutation matrices to inversion words thus
obtained is \EM{injective}, whence it is a bijection.
\end{proof}
Now we simply employ
Striker's bijection \cite{Striker:2011} between inversions words and \DPP s without special parts:
This bijection maps an inversion word $\pas{a_1, a_2,\dots, a_{n-1}}$ to a \DPP\ $\pi$ of dimension $n$
without special parts and with precisely $a_i$ parts $\pas{n+1-i}$, whence we obtain
\bit
\item $\stati\of \pi = a_1+a_2+\cdots+a_{n-1}$,
\item $\statp\of \pi = a_1$,
\item $\statb\of \pi = a_2 + \Iverson{a_2 + a_1 \geq n-1}$\footnote{Here we use Iverson's notation again.}.
\eit
(the last assertion is due to the simple characterization \eqref{eq:path-non-special}).

We claim that our bijection does, in fact, respect the quadruplet of statistics $\quadruplet$.
\begin{proof}
The assertion is trivial for $\stats$ (we consider only objects with $\stats\equiv0$) and obvious for $\stati$ and $\statp$,
so it remains to show this for $\statb$:
Letting $\epsilon\defeq\Iverson{\statb\of A + a_1 \geq n}$, we have
$a_2=\statb\of A -\epsilon$ and thus obtain
$$
\statb\of\pi =
{\statb\of A - \epsilon} + \Iverson{\pas{\statb\of A - \epsilon} + a_1 \geq n-1}.
$$
We have to consider two cases:
\bit
\item $\epsilon=1\iff\statb\of A + a_1 \geq n$: This implies $\Iverson{\statb\of A -\epsilon + a_1\geq n-1}=\epsilon$.
\item $\epsilon=0\iff\statb\of A + a_1 < n$: Observe that $\statb\of A + a_1 = \statb\of A + \statp\of A =n-1$
	ist \EM{not possible} for a \EM{permutation matrix} $A$, since this would imply that the entries $1$ in the first and
	last row appear in the \EM{same} column. So we must have $\statb\of A + a_1 < n-1$ in this case, whence we again
	obtain $\Iverson{\statb\of A -\epsilon + a_1\geq n-1}=\epsilon$.
\eit
So in both cases, we have $\statb\of\pi = \statb\of A$.
\end{proof}

\begin{figure}
\caption{Illustration of the steps of the bijective construction.}
Starting with the $7\times 7$ permutation matrix corresponding to the permutation
$\sigma=\pas{5,2,1,7,3,6,4}$, the pictures show the $6$ steps of the bijective
construction: The rows and columns to be deleted are indicated by thick gray lines.
The inversion word thus obtained is $\pas{4,2,1,1,0,1}$.
\label{fig:bij-example}
\begin{center}
\input ./save-graphics/ASM7step1
\hfil
\input ./save-graphics/ASM6step2
\hfil
\input ./save-graphics/ASM5step3
\end{center}

\begin{center}
\input ./save-graphics/ASM4step4
\hfil
\input ./save-graphics/ASM3step5
\hfil
\input ./save-graphics/ASM2step6
\end{center}
\end{figure}

\secB{An open question}
\label{sec:question}

Our construction was successful because we have the ``translation''
\eqref{eq:path-non-special}
of the (somewhat complicated) condition 
\begin{quote}
\DPP\ $\pi$ of dimension $n$ has a row of length $n-1$
\end{quote}
for \DPP s \EM{without} special parts to an obvious (and quite simple) corresponding
condition for \ASM s \EM{without} entries $-1$.

In the search for a ``natural'' bijection between \ASMes\ and \DPP s, it might be helpful to identify such ``corresponding condition''
for \EM{general} \ASM s.



\bibliographystyle{plain}
\bibliography{paper}

\end{document}